\documentclass[12pt,oneside]{amsart}
\usepackage{amsmath}
\usepackage[mathscr]{euscript}
\usepackage{epsfig}
\usepackage{graphics}
\usepackage{amsbsy}
\usepackage{amsfonts}
\usepackage{amssymb}
\usepackage{graphicx}

 \newcommand{\bd}{\begin{definition}}
 \newcommand{\ed}{\end{definition}}
 \newcommand{\bt}{\begin{theorem}}
 \newcommand{\et}{\end{theorem}}
 \newcommand{\bp}{\begin{proposition}}
 \newcommand{\ep}{\end{proposition}}
 \newcommand{\bl}{\begin{lemma}}
 \newcommand{\el}{\end{lemma}}
 \newcommand{\bpr}{\begin{proof}}
 \newcommand{\epr}{\end{proof}}
 \newcommand{\bc}{\begin{corollary}}
 \newcommand{\ec}{\end{corollary}}
 \newcommand{\be}{\begin{example}}
 \newcommand{\ee}{\end{example}}
\newcommand{\Mm}{\mathcal{M}}

\usepackage{amssymb}
\usepackage{mathrsfs}
\newtheorem{theorem}{{\bf Theorem}}[section]
\newtheorem{lemma}[theorem]{{\bf Lemma}}
\newtheorem{corollary}[theorem]{{\bf Corollary}}
\newtheorem{proposition}[theorem]{{\bf Proposition}}
\newtheorem{definition}[theorem]{{\bf Definition}}
 \makeatletter
      \def\@setcopyright{}
      \def\serieslogo@{}
      \makeatother
\title[$c$-nilpotent multiplier of $n$-Lie algebras]{On the dimension of $c$-nilpotent multiplier of $n$-Lie algebras}
\normalsize{
\author{Farshid Saeedi$^*$ and Seyedeh Nafiseh Akbarossadat}
}

\begin{document}
\maketitle

\noindent
\textbf{Abstract.} 
Let $L$ be a finite-dimensional $n$-Lie algebra with free presentation $F/R$. Then the concept of $c$-nilpotent multiplier of $L$, denoted by $\Mm^{(c)}(L)$, is defined as follows: 
\[\Mm^{(c)}(L)=\dfrac{\gamma_{c+1}(F)\cap R}{\gamma_{c+1}(R,F,\dots,F)}.\] 
In this paper, we obtain some inequalities and certain bounds for the dimension of $\Mm^{(c)}(L)$ by using the basic commutators. Also, we discuss the relationship between the dimension of the $c$-nilpotent multiplier of $L$ and the $c$-nilpotent multiplier of some factor of $L$.     
We further obtain an inequality between dimensions of $c$-nilpotent multiplier of $n$-Lie algebra and non-abelian tensor (exterior) product of a central ideal by its abelianized factor $n$-Lie algebra. Finally, we also determine the dimension and structure of $c$-nilpotent multipliers Heisenberg $n$-Lie algebras, which can be a useful tool for determining the dimension of the multiplier of nilpotent $n$-Lie algebras of class $2$.

\ \\


\section{Introduction and Preparatory}

\indent 
~~~
%
%
%
Ellis \cite{Ellis1991} developed an analogous theory of non-abelian tensor products for Lie algebras (see also \cite{Ellis1987}). Using tensor (exterior) products of Lie algebras, Ellis described the universal central 
extension of Lie algebras. 




In 1986, Filippov \cite{vtf} introduced the notion of $n$-Lie algebras. An \textit{$n$-Lie algebra} over a field $\Bbb F$ is a vector space $L$ over $\Bbb F$ along with an anti-symmetric $n$-linear form $[x_1,\ldots,x_n]$ satisfying the Jacobi identity:
\[[[x_1,\ldots,x_n],y_2,\ldots,y_n]=\sum_{i=1}^n[x_1,\ldots,x_{i-1},[x_i,y_2,\ldots,y_n],x_{i+1},\ldots,x_n].\]
Clearly, $n$-Lie algebras are nothing but ordinary Lie algebras when $n=2$. Studying $n$-Lie algebras is important because of their applications in physics and geometry. So far, several papers have been published about their classification.
The author and Saeedi defined \cite{Akbarossadat-Saeedi1,Akbarossadat-Saeedi2,Akbarossadat-Saeedi-3} the concepts of non-abelian tensor (exterior) product of $n$-Lie algebras and proved some results about them. Also, we introduced some bounds for the dimension of the non-abelian tensor square and Schur multiplier of $n$-Lie algebras. So far, various researches have been done in the field of multipliers (see \cite{Akbarossadat-2-nilpotent,Araskhan,Araskhan-Rismanchi,Niroomand-Parvizi-M^2(L),Salemkar-Aslizadeh, Salemkar-Edalatzadeh-Araskhan,Salemkar-Riyahi}).

First, we will review some of the concepts required in this paper.


 
An $n$-Lie algebra $L$ is called nilpotent of class $s$ (for some positive and integer number $s$), if $L^{s+1}=\gamma_{s+1}(L)=0$ and $L^s\neq0$. This is equivalent to $Z_{s-1}(L)\subsetneq Z_s(L)=L$. Then $s$ is said the nilpotency class of $L$, and we write $cl(L)=s$. Note that the nilpotency property of $n$-Lie algebras is closed under subalgebra, ideal, and homomorphic image, but it is not closed under the extended property.  

\bd \label{wedge}
Let $L$ and $P$ be two $n$-algebras over an arbitrary field $\Bbb F$ such that $L$ and $P$ act on each other by the families of $n$-linear maps $\{f_i\}_{1\leq i\leq n-1}$ and $\{g_i\}_{1\leq i\leq n-1}$, respectively, and act on themselves by their brackets. Then $L\wedge P$ is generated by all symbols 
\[l_1\wedge\cdots\wedge l_i\wedge p_{i+1}\wedge\cdots\wedge p_n,\]
for all $l_i\in L$, $p_i\in P$, and $\alpha\in\Bbb F$, in which the following properties are hold: 
\begin{align}
l_1&\wedge \cdots\wedge \alpha l_i+l'_i\wedge  \cdots\wedge  l_j \wedge  p_{j+1}\wedge  \cdots\wedge  p_n\nonumber \\ 
=&\alpha(l_1\wedge \cdots\wedge   l_i\wedge  \cdots\wedge  l_j \wedge  p_{j+1}\wedge  \cdots\wedge  p_n)\nonumber \\
&+(l_1\wedge \cdots\wedge  l'_i\wedge  \cdots\wedge  l_j \wedge  p_{j+1}\wedge  \cdots\wedge  p_n),\nonumber \\ \nonumber \\
%
l_1&\wedge \cdots\wedge  l_j \wedge  p_{j+1}\wedge \cdots\wedge   \alpha p_{j+i}+p'_{j+i}\wedge  \cdots\wedge  p_n\nonumber \\ 
=&\alpha(l_1\wedge \cdots\wedge  l_j \wedge  p_{j+1}\wedge  \cdots\wedge   p_{j+i}\wedge  \cdots\wedge  p_n)\nonumber \\
&+(l_1\wedge \cdots\wedge  l_j \wedge  p_{j+1}\wedge  \cdots\wedge  p'_{j+i}\wedge  \cdots\wedge  p_n), \nonumber \\ \nonumber \\
%
[l_1&,\ldots,l_n]\wedge  l'_2\wedge \cdots\wedge  l'_i\wedge  p_1\wedge \cdots\wedge  p_{n-i}\nonumber \\ 
=&\sum_{j=1}^n(-1)^{n-j}l_1\wedge \cdots\wedge  l_{j-1}\wedge  l_{j+1}\wedge \cdots\wedge  l_n\wedge  f_i(l_j,l'_2,\ldots,l'_i,p_1,\ldots ,p_{n-i}),\label{wedge3} \\ \nonumber \\ 
 l_1&\wedge \cdots\wedge  l_i\wedge  p_1\wedge \cdots\wedge  p_{n-i-1}\wedge  [p'_1,\ldots,p'_n]\nonumber \\ 
=&\sum_{j=1}^n(-1)^{n-j-2}g_i(p'_j,l_1,\ldots,l_i,p_1,\ldots,p_{n-i-1})\wedge  p'_1\wedge \cdots\wedge  p'_{j-1}\wedge  p'_{j+1}\wedge \cdots\wedge  p'_n, \nonumber\\ \nonumber \\ 
l_1&\wedge\cdots\wedge l_{j-1}\wedge l^*\wedge\dots\wedge l_{k-1}\wedge l^*\wedge\cdots \wedge l_i\wedge p_1\wedge\cdots \wedge p_{n-i}=0,\nonumber\\
l_1&\wedge\cdots \wedge l_i\wedge p_1\wedge\cdots\wedge p_{j-1}\wedge p^*\wedge\cdots\wedge p_{k-1}\wedge p^*\wedge \cdots \wedge p_{n-i}=0.\nonumber
\end{align}
\begin{align*}
{\small \begin{aligned}
\Big[&(l_1^1\wedge \cdots\wedge  l_{i_1}^1\wedge  p_1^1\wedge \cdots\wedge  p_{n-i_1}^1),(l_1^2\wedge \cdots\wedge  l_{i_2}^2\wedge  p_1^2\wedge \cdots\wedge  p_{n-i_2}^2),\ldots,\\ 
&(l_1^{n-1}\wedge \cdots\wedge  l_{i_{n-1}}^{n-1}\wedge  p_1^{n-1}\wedge \cdots\wedge  p_{n-i_{n-1}}^{n-1}),(l_1^n\wedge \cdots\wedge  l_{i_n}^n\wedge  p_1^n\wedge \cdots\wedge  p_{n-i_{n}}^n)\Big]\\
&=\frac{1}{2^{n-1}}\Big((-1)^{\sum_{\substack{s=1}}^{n-1}i_s(n-i_s)}\big(g_{n-i_1}(p_1^1,\ldots,p_{n-i_1}^1,l_1^1,\ldots,l_{i_1}^1)\\
&\hspace{1.7cm}\wedge  f_{i_2}(l_1^2,\ldots,l_{i_2}^2,p_1^2,\ldots,p_{n-i_2}^2)\\
&\hspace{1.7cm}\wedge \cdots\wedge  f_{i_{n-1}}(l_1^{n-1},\ldots,l_{i_{n-1}}^{n-1},p_1^{n-1},\ldots,p_{n-i_{n-1}}^{n-1})\\
&\hspace{1.7cm}\wedge  f_{i_n}(l_1^n,\ldots ,l_{i_n}^n,p_1^n,\ldots,p_{n-i_{n}}^n)\big)\\ 
&\quad+\sum_{r=2}^{n-1}(-1)^{\sum_{\substack{s=1\\s\neq r}}^{n-1}i_s(n-i_s)}\big(g_{n-i_1}(p_1^1,\ldots,p_{n-i_1}^1,l_1^1,\ldots,l_{i_1}^1)\\
&\hspace{1cm}\wedge   g_{n-i_r}(p_1^r,\ldots,p_{n-i_1}^r,l_1^r,\ldots,l_{i_1}^r)\wedge  f_{i_2}(l_1^2,\ldots,l_{i_2}^2,p_1^2,\ldots,p_{n-i_2}^2) \\
&\hspace{1cm}\wedge \cdots\wedge  f_{i_{r-1}}(l_1^{r-1},\ldots,l_{i_2}^{r-1},p_1^{r-1},\ldots,p_{n-i_2}^{r-1})\\
&\hspace{1cm}\wedge  f_{i_{r+1}}(l_1^{r+1},\ldots,l_{i_2}^{r+1},p_1^{r+1},\ldots,p_{n-i_2}^{r+1})\\ 
&\hspace{1cm}\wedge \cdots\wedge  f_{i_{n-1}}(l_1^{n-1},\ldots,l_{i_2}^{n-1},p_1^{n-1},\ldots,p_{n-i_2}^{n-1})\\
&\hspace{1cm}\wedge f_{i_n}(l_1^n,\ldots,l_{i_2}^n,p_1^n,\ldots,p_{n-i_2}^n)\big)\\
&\quad+\cdots+(-1)^{i_1(n-i_1)}\big(g_{n-i_1}(p_1^1,\ldots,p_{n-i_1}^1,l_1^1,\ldots,l_{i_1}^1)\\
&\hspace{1cm}\wedge  g_{n-i_2}(p_1^2,\ldots,p_{n-i_2}^2,l_1^2,\ldots,l_{i_2}^2)\\
&\hspace{1cm}\wedge \cdots\wedge  g_{n-i_{n-1}}(p_1^{n-1},\ldots,p_{n-i_{n-1}}^{n-1},l_1^{n-1},\ldots,l_{i_{n-1}}^{n-1})\\
&\hspace{1cm}\wedge  g_{n-i_{n}}(p_1^{n},\ldots ,p_{n-i_{n}}^{n},l_1^{n},\ldots,l_{i_{n}}^{n})\big)\Big),
\end{aligned}}
\end{align*}
\ed
In what  follows, we recall the definition of crossed modules of $n$-Lie algebras. 
\bd[Crossed module] \label{crossed module}
A crossed module is a homomorphism of Lie  $n$-algebras $\mu : L\longrightarrow P$ 
together with an action $\{g_i\}_{1\leq i\leq n-1}$ of $P$ on $L$  satisfying the following conditions:
\begin{enumerate}
\item 
$\mu$  is compatible with the action of $P$ on $L$, that is, 
\[\mu\left(g_i(p_1,\ldots, p_i,l_1,\ldots,l_{n-i})\right)=[p_1,\ldots,p_{n-i},\mu(l_1),\ldots, \mu(l_i)],\]
for all $1\leq i\leq n-1$.
\item 
For all $1\leq i\leq n-1$, 
\[g_i(\mu(l_1),\ldots,\mu(l_i),l_{i+1},\ldots,l_{n})=[l_1,\ldots,l_i,l_{i+1},\ldots,l_n].\]
\item 
For all $1\leq i\leq n-1$ and $i+1\leq j\leq n$,
\[g_i(p_1,\ldots,p_i,l_{i+1},\ldots,l_{n})=(-1)^{j-i-1}g_{i+1}(p_1,\ldots,p_i,\mu(l_j),l_{i+1},\ldots,l_{j-1},l_{j+1},\ldots,l_{n}).\]
\end{enumerate}
\ed 
It is easy to check that the kernel of $\mu$  is in the center of $L$ and that its image is an ideal in $P$.
Furthermore, the $n$-Lie algebra  $\mathrm{Im}(\mu)$ acts trivially on the  center $Z(L)$, and so trivially on 
$\ker\mu$. Hence $\ker \mu$ inherits an action of $P/\mathrm{Im}(\mu)$ making $\ker\mu$ a  representation of the $n$-Lie algebra $P/\mathrm{Im}(\mu)$.

For the first time, the authors \cite{Akbarossadat-Saeedi-3} introduced the concept of free $n$-Lie algebras and, then in \cite{Akbarossadat-Basic}, defined the concept of basic commutators of weight $w$ in $d$-dimensional $n$-Lie algebras. Also, we proved some of its properties and the formula to calculate the number of them. In what follows, we review it. 
\begin{theorem}[\cite{Akbarossadat-Basic}]\label{main-formula-weightw}
Let the set $X=\{x_i|x_{i+1}>x_i;~i=1,2,\dots,d\}$ be an ordered set and a basis for the free $n$-Lie algebra $F$ 
 and let 		$w$ be a positive integer number. Then the number of basic commutators of weight $w$ is 
\begin{equation}\label{basicweightw}
l_d^n(w)=\sum_{j=1}^{\alpha_0}\beta_{j^*}\left(\sum_{i=2}^{w-1}\alpha_i{{{{d}\choose{n-1}}}\choose{w-i}}\right),
\end{equation}
where $\alpha_0={{d-1}\choose{n-1}}$, $\alpha_i$, ($2\leq i\leq w-1$) is the coefficient of the $(i-2)$th sentence in Newton's binomial expansion $(a+b)^{w-3}$, $({\text{i.e.}}~ \alpha_i={{w-3}\choose{i-2}})$ and if ${{k-1}\choose{n-1}}+1\leq j\leq {{k}\choose{n-1}}$ (for $k=n-1,n,n+1,n+2,\dots,d-1$), then $j^*={{k-1}\choose{n-1}}+1$ and $\beta_{j^*}=(d-n-j^*+2)$.
\end{theorem}
The following theorem expresses one of the most important and main applications of basic commutators of $n$-Lie algebras. 
\bt[\cite{Akbarossadat-Basic}] \label{F^i/F^j}
Let $F$ be a free $n$-Lie algebra and let $F^i$ be the $i$th term of the lower central series of $F$, for each $i\in\Bbb N$. Then $\dfrac{F^i}{F^{i+c}}$ is abelian of dimension $\sum\limits_{j=0}^{c-1}l_d^n(i+j)$, where $c=1,2,\ldots$.   
\et 
%
The next proposition determines the dimension of $c$-nilpotent multiplier of all abelian $n$-Lie algebras. 
\bp[\cite{Akbarossadat-2-nilpotent}] \label{pro1.2} 
Let $L$ be an abelian $n$-Lie algebra with  dimension finite $d$. Then 
$\dim\Mm^{(c)}(L) = l^n_d(c + 1)$.
 In particular,
$\dim\Mm(L) =l_d^n(2)=\dfrac{1}{2}d(d-1)$.
\ep 
\section{Preliminary results}

\indent 

In this section, we try to present some important and practical tools in order to prove the main results of this paper. 

\bl \label{lem-0}
Let $0\longrightarrow R\longrightarrow F\longrightarrow L\longrightarrow0$ be a free presentation of the $n$-Lie algebra $L$ and let $M$ be an ideal of $L$ with free presentation $S/R$, for some ideal $S$ of $F$. Then 
\[\gamma_{c+1}(L)\cap M\cong\dfrac{(\gamma_{c+1}(F)\cap S)/\gamma_{c+1}(R,F,\dots,F)}{\Mm^{(c)}(L)}.\]
\el
\bpr 
By the definition of the lower central series, we know that 
\begin{align*}
\gamma_{c+1}(L)\cap M= \gamma_{c+1}(F/R)\cap (S/R)
&=\dfrac{\gamma_{c+1}(F)+R}{R}\cap \dfrac{S}{R}\\
&=\dfrac{(\gamma_{c+1}(F)+R)\cap S}{R}\\
&=\dfrac{(\gamma_{c+1}(F)\cap S)+(R\cap S)}{R}\\
&=\dfrac{(\gamma_{c+1}(F)\cap S)+R}{R}\\
&\cong\dfrac{(\gamma_{c+1}(F)\cap S)}{(\gamma_{c+1}(F)\cap S)\cap R}\\
&=\dfrac{(\gamma_{c+1}(F)\cap S)}{\gamma_{c+1}(F)\cap R}\\
&\cong\dfrac{(\gamma_{c+1}(F)\cap S)/\gamma_{c+1}(R,F,\dots,F)}{(\gamma_{c+1}(F)\cap R)/\gamma_{c+1}(R,F,\dots,F)}\\
&=\dfrac{(\gamma_{c+1}(F)\cap S)/\gamma_{c+1}(R,F,\dots,F)}{\Mm^{(c)}(L)}.
\end{align*}
\epr 

\bl \label{lem2.3(2)}
Let $L$ be an $n$-Lie algebra with a free presentation $L\cong F/R$, and let $M$ be its central ideal and $M\cong S/R$, for some ideals $S$ of the free $n$-Lie algebra $F$. Then $\gamma_{c+1}(S,F,\dots,F)/(\gamma_{c+1}(R,F,\dots,F)+\gamma_{c+1}(S))$ is a homomorphic image of 
\[(L/M)^{ab}\ ^c\otimes M=\underbrace{(L/M)^{ab}\otimes\dots\otimes (L/M)^{ab}}_{c-times}\otimes M,\]
 where $(L/M)^{ab}=\dfrac{L/M}{(L/M)^2}$ and the tensor product is in sense of tensor product of $n$-Lie algebras.  
\el 
\bpr 
Since $F/R$ is a free presentation of $L$, so there exists an epimorphism $\pi:F\longrightarrow L$ such that $\ker\pi=R$. 
Consider the following map: 
\[\theta:\begin{cases}
\underbrace{(L/M)^{ab}\times\dots\times (L/M)^{ab}}_{(c-1)(n-1)-times}\times M\longrightarrow \gamma_{c+1}(S,F,\dots,F)/(\gamma_{c+1}(R,F,\dots,F)+\gamma_{c+1}(S))\\
(\bar{l}_1,\dots,\bar{l}_n,\bar{l}_{n+1},\dots,\bar{l}_{2n-1},\dots,\bar{l}_{(c-1)(n-1)})\longmapsto\\ \\
 \overline{[\dots[[\bar{f}_1,\dots,\bar{f}_{n-1},m],\bar{f}_{n},\bar{f}_{n+1},\dots,\bar{f}_{2n-2}],\dots],\bar{f}_{(c-2)(n-1)+1},\dots,\bar{f}_{(c-1)(n-1)}],\dots]},
\end{cases}\]
where $\pi^{-1}(l_i)=\bar{f}_i$, for all $i$. 
It is easy to check that $\theta$ is linear, and hence it induces the following epimorphism:
\[(L/M)^{ab}\ ^c\otimes M=\underbrace{(L/M)^{ab}\otimes...\otimes (L/M)^{ab}}_{c-times}\otimes M\rightarrow \gamma_{c+1}(S,F,...,F)/(\gamma_{c+1}(R,F,...,F)+\gamma_{c+1}(S)).\]
\epr

The authors \cite{Akbarossadat-Saeedi1} proved a result similar to the following lemma for non-abelian tensor product of $n$-Lie algebras. By a similar method and reasoning, the following lemma can be proved for the non-abelian exterior (wedge) product of $n$-Lie algebras as well.
\bl \label{lemma1.3}
Let $L$, $P$, $K$, and $Q$ be $n$-Lie algebras such that $L$ and $K$ act compatibly on each other by $\{f_i\}_{1\leq i\leq n-1}$ and $\{g_i\}_{1\leq i\leq n-1}$, and so they act $P$ and $Q$ by $\{h_i\}_{1\leq i\leq n-1}$ and $\{s_i\}_{1\leq i\leq n-1}$, respectively. Moreover, suppose that $\sigma_1:L\longrightarrow P$ and $\sigma_2:K\longrightarrow Q$ are $n$-Lie homomorphisms preserving the actions in the sense that
\begin{align*}
&\sigma_1(g_i(a_1,\ldots ,a_n))=s_i(\sigma_{r}(a_1),\ldots ,\sigma_{r}(a_n)),\\
&\sigma_2(f_i(a_1,\ldots ,a_n))=h_i(\sigma_{r}(a_1),\ldots ,\sigma_{r}(a_n)),
\end{align*}
for all $a_1,\ldots ,a_n\in L\cup K$, where $r=1$ if $a_i\in L$, and $r=2$ if $a_i\in K$. 
Then there is a unique homomorphism (up to sign) 
{\small{\[\begin{array}{rcl}
\sigma_1\wedge \sigma_2:L\wedge K&\longrightarrow &P\wedge Q\\
\sum\limits_{i=1}^{n-1}l_1^i\wedge \cdots\wedge l_i^i\wedge k_{i+1}^i\wedge \cdots\wedge k_n^i&\longmapsto &\sum\limits_{i=1}^{n-1}\sigma_1(l_1^i)\wedge \cdots\wedge \sigma_1(l_i^i)\wedge \sigma_2(k_{i+1}^i)\wedge \cdots\wedge \sigma_2(k_n^i).
\end{array}\]}}
\el 
\bpr 
The proof is  similar to  \cite[Proposition 4.5]{Akbarossadat-Saeedi1}. 
\epr 
\bl \label{lemma01}
Let $F$ be a free $n$-Lie algebra and let $S$ and $I$ be two ideals of $F$. Then 
\[[I,\gamma_j(S,F,\dots,F),F,\dots,F]\subseteq\gamma_j(I,F,\dots,F),\]
for all $j\geq 1$. 
\el 
\bpr 
We prove this by induction on $j$. The result is trivially true for $j=1$. Assume that the result holds for $j\geq1$. Then for any ideal $I$ in $F$, by the definition of $\gamma_{j+1}$ and the Jacobi identity, we have 
\begin{align*}
[I,\gamma_{j+1}(S,F,\dots,F),F,\dots,F]&=[\gamma_{j+1}(S,F,\dots,F),I,F,\dots,F]\\
&=[[\gamma_{j}(S,F,\dots,F),\underbrace{F,\dots, F}_{n-1}],I,\underbrace{F,\dots,F}_{n-2}]\\
&=[[\gamma_{j}(S,F,\dots,F),I,\underbrace{F,\dots,F}_{n-2}],\underbrace{F,\dots, F}_{n-1}]\\
&\quad+[[\gamma_{j}(S,F,\dots,F),[F,I,\underbrace{F,\dots,F}_{n-2}],\underbrace{F,\dots, F}_{n-2}]\\
&\quad+\dots\\
&\quad+[\gamma_{j}(S,F,\dots,F),\underbrace{F,\dots,F}_{n-2}, [F,I,\underbrace{F,\dots,F}_{n-2}]]\\
&\subseteq [\gamma_{j}(I,F,\dots,F),F,\dots,F]+[\gamma_{j}(S,F,\dots,F),I,F,\dots,F]\\
&\subseteq \gamma_{j}([I,F,\dots,F],F,\dots,F)+\gamma_{j+1}(I,F,\dots,F) \\
&=\gamma_{j+1}(I,F,\dots,F).
\end{align*}
\epr 
Let $L$ be an $n$-Lie algebra with a free presentation $F/R$, for some free $n$-Lie algebra $F$ and an ideal $R$ of $F$ with the $n$-Lie algebra epimorphism $\pi$, and let $M$ be an ideal of $L$ with a free presentation $S/R$, for some ideal $S$ of $F$. Since $R\subseteq S\subseteq F$, so $\gamma_j(R,F,\dots,F)\subseteq \gamma_j(S,F,\dots,F)$, and hence we can put $T_j=\dfrac{\gamma_j(S,F,\dots,F)}{\gamma_j(R,F,\dots,F)}$, for all $j\geq 1$. The following lemma shows that $T_j$ and $L$ act on each other compatibly. 
\bl \label{lemma02}
By the above assumptions and notations, $T_j$ and $L$ act compatibly on each other.
\el 
\bpr 
Define the two family of maps $\{g_i\}$ and $\{h_i\}$ as follows:
\[\begin{cases} 
g_i:L^{\times i}\times T_j^{\times n-i}\longrightarrow T_j\\
g_i\left(l_1,\dots,l_i,\bar{t_{i+1}},\dots,\bar{t_n}\right)=[f_1,\dots,f_i,t_{i+1},\dots,t_n]+\gamma_j(R,F,\dots,F)
\end{cases}\]
and 
\[\begin{cases} 
h_i:T_j^{\times i}\times L^{\times n-i}\longrightarrow L\\
h_i\left(\bar{t_{1}},\dots,\bar{t_i},l_{i+1},\dots,l_n\right)=\pi([t_1,\dots,t_i,l_{i+1},\dots,l_n])
\end{cases}\]  
where $l_p\in L$, $x_q\in \gamma_j(S,F,\dots,F)$, and $f_p=\pi^{-1}(l_p)$, for all $p$ and $q$. By Lemma \ref{lemma01}, it is easy to check that $T_j$ and $L$ act compatibly on each other.
\epr 

The following proposition is useful in  some of our  results and future investigation.
\bp \label{pro1.4}
Let $L$ be an $n$-Lie algebra and let $M$ be its ideal. Then 
\begin{enumerate}
\item 
$T_{j+1}$ is an isomorphic image of $M\wedge^j L$. 
\item 
if $M$ is an $r$-central ideal of $L$, then $T_{j+1}$ is an isomorphic image of $M\wedge^jK$, where $K=\dfrac{L}{\gamma_{r+1}(L)}$.  
\end{enumerate}
\ep 
\bpr 
\begin{enumerate}
\item 
Consider the map $\lambda:T_j\longrightarrow L$ by $\lambda(\bar{t})=\pi(t)$. Since $\pi$ is a well-defined $n$-Lie homomorphism, so $\lambda$ is too. The $n$-Lie homomorphism $\lambda$ together with the defined action of $L$ on $T_j$ in Lemma \ref{lemma02} and also, the identity map $id_L:L\longrightarrow L$ are crossed modules of $n$-Lie algebras. 
The family of bilinear maps $\phi_j^i:T_j^{\times i}\times L^{\times n-i}\stackrel{epi.}{\longrightarrow}T_{j+1}$ given by $\phi_(\bar{x}_1,\dots,\bar{x}_i,l_{i+1},\dots,l_{n}) = [x_1,\dots,x_i,f_{i+1},\dots,f_n]+\gamma_{j+1}(R,F,\dots,F)$, is an $n$-multiplying. By   Theorem 3.5 of \cite{Akbarossadat-Saeedi1} consequently, this family induces an epimorphism $\Phi_j:T_j \wedge L\longrightarrow T_{j+1}$. Thus, by using Lemma \ref{lemma1.3}, we obtain an epimorphism $\Phi_j\wedge id_L:(T_j\wedge L)\wedge L\longrightarrow T_{j+1}\wedge L$. Now, we continue by induction on $j$.

\item 
By the assumption, $\gamma_{n+1}(S,F,\dots,F)\subseteq R$ and then 
\[[\gamma_j(S,F,\dots,F),\gamma_{n+1}(F),F,\dots,F]\subseteq \gamma_{j+n+1}(S,F,\dots,F)\subseteq \gamma_{j+1}(R,F,\dots,F).\]
So, the action of $\gamma_{n+1}(F)$ on $T_j$ is trivial and the $n$-Lie algebras $L/\gamma_{n+1}(L)$ and $T_j$ act compatibly on each other by the induced actions $\bar{g}_i$'s and $\bar{h}_i$'s of defined $g_i$'s and $h_i$'s in Lemma \ref{lemma02}; that is,
\[\begin{cases} 
\bar{g}_i:L/\gamma_{n+1}(L)^{\times i}\times T_j^{\times n-i}\longrightarrow T_j\\
\bar{g}_i\left(\bar{l}_1,\dots,\bar{l}_i,\bar{t}_{i+1},\dots,\bar{t}_n\right)=g_i\left(l_1,\dots,l_i,\bar{t}_{i+1},\dots,\bar{t}_n\right)\\
\qquad\qquad\qquad\qquad\qquad\quad=[f_1,\dots,f_i,t_{i+1},\dots,t_n]+\gamma_j(R,F,\dots,F)
\end{cases}\]
and 
\[\begin{cases} 
\bar{h}_i:T_j^{\times i}\times L/\gamma_{n+1}(L)^{\times n-i}\longrightarrow L/\gamma_{n+1}(L)\\
\bar{h}_i\left(\bar{t_{1}},\dots,\bar{t_i},\bar{l}_{i+1},\dots,\bar{l}_n\right)=\pi([t_1,\dots,t_i,l_{i+1},\dots,l_n])
\end{cases}\]  
where $l_p\in L$, $x_q\in \gamma_j(S,F,\dots,F)$, and $f_p=\pi^{-1}(l_p)$, for all $p$ and $q$. 

Now, as part (1), there exist epimorphisms 
$T_{j}\wedge K\longrightarrow T_{j+1}$ and $M\wedge^j K\longrightarrow T_{j+1}$ and hence the proof is completed. 
\end{enumerate}
\epr 
\bc \label{cor1.5}
By the assumptions and notations in Proposition \ref{pro1.4}, there exists an epimorphism $\kappa$ from $\ker\mu^j_M$ onto $\dfrac{R\cap \gamma_{j+1}(S,F,\dots,F)}{\gamma_{j+1}(R,F,\dots,F)}$.
\ec 
\bpr 
The proof is similar to   \cite[Corollary 1.5]{Salemkar-Edalatzadeh-Araskhan}. 
\epr 
\section{On the $c$-nilpotent multiplier of $n$-Lie algebras} 
~~~In this section, first we introduce some inequalities for the dimension of the $c$-nilpotent multiplier of the $n$-Lie algebras
and their factors. Then, we show that  every $c$-perfect $n$-Lie algebra has at least one cover. Also, we give some upper bounds for the dimension of $n$-nilpotent multiplier of $n$-Lie algebras using the basic commutators. 

The following lemma is one of the important tools for obtaining the required inequalities, which is similar to the work of Araskhan and Rismanchian  \cite{Araskhan-Rismanchi} for the Lie algebras  and the work of Jones \cite{[10]-2} for the group case.

\bl \label{lem2.1(2)}
Let $L$ be an $n$-Lie algebra with dimension finite, and let $M$ be its ideal. Then there exists an $n$-Lie algebra $K$ with its ideal $S$ such that 
\begin{enumerate}
\item 
$\gamma_{c+1}(L)\cap M=K/S$;
\item 
$S\cong\Mm^{(c)}(L)$;
\item 
$\Mm^{(c)}(L/M)$ is a homomorphic image of $K$.
\item 
If, in addition, $M$ is a $c$-central ideal of $L$, then $\gamma_{c+1}(L)\cap M$ is a homomorphic image
of $\Mm^{(c)}(L/M)$.
\end{enumerate}
\el 
\bpr 
\begin{enumerate}
\item 
Since $\gamma_{c+1}(F)\cap R\unlhd \gamma_{c+1}(F)\cap S$, so by  Lemma \ref{lem-0}, it is enough to put $K/M=\dfrac{(\gamma_{c+1}(F)\cap S)/\gamma_{c+1}(R,F,\dots,F)}{(\gamma_{c+1}(F)\cap R)/\gamma_{c+1}(R,F,\dots,F)}$. 
\item 
By the definition of the $c$-nilpotent multiplier of $L$ and Lemma \ref{lem-0}, we obtain  $S=\Mm^{(c)}(L)$. 
\item 
Consider the free presentation
\[L/M=\dfrac{F/R}{S/R}\cong F/S,\]
and hence 
\[\Mm^{(c)}(L/M)=\dfrac{\gamma_{c+1}(F)\cap S}{\gamma_{c+1}(S,F,\dots,F)}.\]
Now, consider the map $\alpha:(\gamma_{c+1}(F)\cap S)/\gamma_{c+1}(R,F,\dots,F)\longrightarrow \dfrac{\gamma_{c+1}(F)\cap S}{\gamma_{c+1}(S,F,\dots,F)}$. 
It is easy to check that $\alpha$ is an epimorphism of $n$-Lie algebras with the kernel ideal $\gamma_{c+1}(S,F,\dots,F)/\gamma_{c+1}(R,F,\dots,F)$. 
\item 
Since $M$ is a $c$-central ideal, so $\gamma_{c+1}(M,L,\dots,L)=1$. Hence $\dfrac{\gamma_{c+1}(S,F,\dots,F)+R}{R}=R$. Thus $\gamma_{c+1}(S,F,\dots,F)\subseteq R$. Therefore, $\gamma_{c+1}(S,F,\dots,F)\subseteq \gamma_{c+1}(F)\cap R$ and $\gamma_{c+1}(S,F,\dots,F)\subseteq \gamma_{c+1}(F)\cap S$.   
Also, by the proof of Lemma \ref{lem-0}, we know that 
\begin{align*}
\gamma_{c+1}(L)\cap M\cong\dfrac{\gamma_{c+1}(F)\cap S}{\gamma_{c+1}(F)\cap R}\cong\dfrac{\dfrac{\gamma_{c+1}(F)\cap S}{\gamma_{c+1}(S,F,\dots,F)}}{\dfrac{\gamma_{c+1}(F)\cap R}{\gamma_{c+1}(S,F,\dots,F)}}=\dfrac{\Mm^{(c)}(L/M)}{\dfrac{\gamma_{c+1}(F)\cap R}{\gamma_{c+1}(S,F,\dots,F)}}.
\end{align*}
Thus the proof is completed by considering the epimorphism $\alpha:\Mm^{(c)}(L/M)\longrightarrow \gamma_{c+1}(L)\cap M=\dfrac{\Mm^{(c)}(L/M)}{\dfrac{\gamma_{c+1}(F)\cap R}{\gamma_{c+1}(S,F,\dots,F)}}$. 
\end{enumerate}
\epr 

The following corollary follows immediately from the previous lemma.
\bc \label{cor2.2(2)} 
Let $L$ be an $n$-Lie algebra with dimension finite, and let $M$ be  its ideal. Then we have the following inequality:
\[\dim\Mm^{(c)}(L/M)\leq \dim\Mm^{(c)}(L)+\dim(\gamma_{c+1}(L)\cap M).\]
\ec 
\bpr 
From Lemma \ref{lem2.1(2)}, we have $\gamma_{c+1}(L)\cap M=K/S$ and $S\cong\Mm^{(c)}(L)$. Thus 
\[\dim(\gamma_{c+1}(L)\cap M)=\dim K-\dim S, \qquad \dim S=\dim\Mm^{(c)}(L).\]
Hence, 
\[\dim(\gamma_{c+1}(L)\cap M)=\dim K-\dim\Mm^{(c)}(L).\]
On the other hand, since $\Mm^{(c)}(L/M)$ is a homomorphic image of $K$, so $\dim\Mm^{(c)}(L/M)\leq \dim K$. Therefore, 
\[\dim\Mm^{(c)}(L/M)\leq \dim\Mm^{(c)}(L)+\dim(\gamma_{c+1}(L)\cap M).\]
\epr 

\bt \label{the2.4(2)} 
Let $L$ be an $n$-Lie algebra with dimension finite and free presentation $F/R$, and let $M$ be a central ideal of $L$ with free presentation $S/R$, for some ideal $S$ of $F$. Then
\[\dim\Mm^{(c)}(L)+\dim(\gamma_{c+1}(L)\cap M)\leq\dim\Mm^{(c)}(L/M)+\dim\Mm^{(c)}(M)+\dim((L/M)^{ab}\ ^c\otimes M).\]
\et 
\bpr 
Since $M$ is central, so $\gamma_2(S,F,\dots,F)\subseteq R$. By Corollary \ref{cor2.2(2)} and the proof of  part (3) of Lemma \ref{lem2.1(2)}, we know that $\Mm^{(c)}(L/M)$ is an isomorphic image of $(\gamma_{c+1}(F)\cap S)/\gamma_{c+1}(R,F,\dots,F)$ with the kernel ideal $\gamma_{c+1}(S,F,\dots,F)/\gamma_{c+1}(R,F,\dots,F)$. Also, we have 
\[\dim\Mm^{(c)}(L/M)+\dim\dfrac{\gamma_{c+1}(S,F,\dots,F)}{\gamma_{c+1}(R,F,\dots,F)}= \dim\Mm^{(c)}(L)+\dim(\gamma_{c+1}(L)\cap M).\]
On the other hand, we know that 
\[\dfrac{\dfrac{\gamma_{c+1}(S,F,\dots,F)}{\gamma_{c+1}(R,F,\dots,F)}}{\dfrac{\gamma_{c+1}(R,F,\dots,F)+\gamma_{c+1}(S)}{\gamma_{c+1}(R,F,\dots,F)}}=\dfrac{\gamma_{c+1}(S,F,\dots,F)}{\gamma_{c+1}(R,F,\dots,F)+\gamma_{c+1}(S)}.\]
Hence 
\[\dim \dfrac{\gamma_{c+1}(S,F,\dots,F)}{\gamma_{c+1}(R,F,\dots,F)}-\dim\dfrac{\gamma_{c+1}(R,F,\dots,F)+\gamma_{c+1}(S)}{\gamma_{c+1}(R,F,\dots,F)}=\dim\dfrac{\gamma_{c+1}(S,F,\dots,F)}{\gamma_{c+1}(R,F,\dots,F)+\gamma_{c+1}(S)}.\]
Also, by this fact that $M$ is a central ideal,  we have 
\begin{align}
&\gamma_{c+1}(S)\subseteq \gamma_2(S,F,\dots,F)\subseteq R\label{equation(1)}\\ 
&\gamma_2(R,F,\dots,F)\subseteq R.\label{equation(2)}
\end{align}
Therefore, it follows from equations \eqref{equation(1)} and \eqref{equation(2)}  that 
\begin{equation}\label{equation(3)}
\gamma_2(R,F,\dots,F)+\gamma_{c+1}(S)\subseteq R.
\end{equation} 
On the other hand, since $R\subseteq S$, so
\begin{equation}\label{equation(4)}
\gamma_2(R,F,\dots,F)\subseteq \gamma_{c+1}(S),
\end{equation}
 and hence equations \eqref{equation(3)} and \eqref{equation(4)} imply that 
\begin{equation}\label{equation(5)}
\gamma_2(R,F,\dots,F)+\gamma_{c+1}(S)\subseteq \gamma_{c+1}(S)\cap R.
\end{equation}
Thus 
\begin{align*}
\dim\Mm^{(c)}(L)+\dim(\gamma_{c+1}(L)\cap M)&=\dim\Mm^{(c)}(L/M)+\dim\dfrac{\gamma_{c+1}(S,F,\dots,F)}{\gamma_{c+1}(R,F,\dots,F)}\\
&=\dim\Mm^{(c)}(L/M)+\dim\dfrac{\gamma_{c+1}(S,F,\dots,F)}{\gamma_{c+1}(R,F,\dots,F)+\gamma_{c+1}(S)}\\
&\quad+\dim \dfrac{\gamma_{c+1}(R,F,\dots,F)+\gamma_{c+1}(S)}{\gamma_{c+1}(R,F,\dots,F)}\\
(\text{by Lemma \ref{lem2.3(2)}})\qquad\qquad &\leq\dim\Mm^{(c)}(L/M)+\dim\left((L/M)^{ab}\ ^c\otimes M\right) \\
&\quad+\dim \dfrac{\gamma_{c+1}(R,F,\dots,F)+\gamma_{c+1}(S)}{\gamma_{c+1}(R,F,\dots,F)}\\
(\text{by \eqref{equation(5)}})\qquad\qquad\qquad \  &\leq \dim\dfrac{\gamma_{c+1}(S)\cap R}{\gamma_{c+1}(R,F,\dots,F)}\\
&=\dim\Mm^{(c)}(L/M)+\dim\Mm^{(c)}(M)\\
&\quad+\dim\left((L/M)^{ab}\ ^c\otimes M\right).
\end{align*}
\epr 

The following proposition plays a fundamental and important role in obtaining  some main results of this section.
\bp \label{pro2.1}
Let $L$ be an $n$-Lie algebra and let $M$ be  its ideal. Then we have the following exact sequences:
\begin{align}
& \ker(\mu_M^c)\longrightarrow\Mm^{(c)}(L)\longrightarrow\Mm^{(c)}(L/M)\longrightarrow\dfrac{M\cap\gamma_{c+1}(L)}{\gamma_{c+1}(M,L,\dots,L)}\longrightarrow 0; \label{sequ.(a)}\\
&M\wedge^c\dfrac{L}{\gamma_{r+1}(L)}\longrightarrow\Mm^{(c)}(L)\longrightarrow\Mm^{(c)}(L/M)\longrightarrow M\cap \gamma_{c+1}(L)\longrightarrow 0,\label{sequ.(b)}
\end{align}
with the condition that $M$ is $r$-central and $r\leq c$. 
\ep 
\bpr 
Let $0\longrightarrow R\longrightarrow F\longrightarrow L\longrightarrow 0$ be a free presentation of $L$ and let $M\cong S/R$, for some ideal $S$ of $F$. Then we have the following exact sequence:
\begin{align}
0\longrightarrow \dfrac{R\cap\gamma_{c+1}(S,F,\dots,F)}{\gamma_{c+1}(R,F,\dots,F)}\longrightarrow \dfrac{R\cap \gamma_{c+1}(F)}{\gamma_{c+1}(R,F,\dots,F)}&\longrightarrow\dfrac{S\cap\gamma_{c+1}(F)}{\gamma_{c+1}(S,F,\dots,F)} \nonumber \\
&\longrightarrow\dfrac{(S+\gamma_{c+1}(F))+R}{\gamma_{c+1}(S,F,\dots,F)+R}\longrightarrow 0. \label{Short-exact-seq.1}
\end{align}
Since $r\leq c$ and $M$ is $c$-central, so 
\begin{equation}\label{equ*}
\gamma_{c+1}(S,F,\dots,F)\subseteq \gamma_{r+1}(S,F,\dots,F)\subseteq R=0_{S/R}.
\end{equation}
We know that $L\cong F/R$ and $M\cong S/R$, and hence $F/S$ is a free presentation of $L/M$.
Thus according to the definition of $c$-nilpotent multiplier of $L/M$, we have 
\begin{align}
\Mm^{(c)}(L)=\dfrac{R\cap\gamma_{c+1}(F)}{\gamma_{c+1}(R,F,\dots, F)},\qquad \Mm^{(c)}(L/M)=\dfrac{S\cap\gamma_{c+1}(F)}{\gamma_{c+1}(S,F,\dots, F)}.\label{equ3}
\end{align}
Also, 
\begin{align}
\dfrac{M\cap \gamma_{c+1}(L)}{\gamma_{c+1}(M,L,\dots, L)}&=\dfrac{(S/R)\cap\gamma_{c+1}(F/R)}{\gamma_{c+1}(S/R,F/R,\dots,F/R)}\nonumber\\ 
&=\dfrac{\dfrac{S+R}{R} \cap\dfrac{\gamma_{c+1}(F)+R}{R}}{\dfrac{\gamma_{c+1}(S,F,\dots,F)+R}{R}}\nonumber\\ 
&=\dfrac{\dfrac{(S\cap\gamma_{c+1}(F))+R}{R}}{\dfrac{\gamma_{c+1}(S,F,\dots,F)+R}{R}}\nonumber\\ 
&=\dfrac{(S\cap\gamma_{c+1}(F))+R}{\gamma_{c+1}(S,F,\dots,F)+R}.\label{equ4}
\end{align} 
By Corollary \ref{cor1.5}, there exists an epimorphism $\psi:\ker(\mu_M^c)\longrightarrow\dfrac{R\cap\gamma_{c+1}(S,F,\dots,F)}{\gamma_{c+1}(R,F,\dots,F)}$, and so 
\begin{equation}\label{equ5}
\dfrac{R\cap\gamma_{c+1}(S,F,\dots,F)}{\gamma_{c+1}(R,F,\dots,F)}\cong \dfrac{\ker(\mu_M^c)}{\ker\psi}.
\end{equation} 
Therefore, by replacing the terms of sequence \eqref{Short-exact-seq.1} by \eqref{equ*}, \eqref{equ3}, \eqref{equ4}, and \eqref{equ5}, the sequence \eqref{sequ.(a)} is obtained. 

On the other hand, from \eqref{equ*}, we have
\begin{equation}\label{equ6}
\dfrac{\gamma_{c+1}(S,F,\dots,F)}{\gamma_{c+1}(R,F,\dots,F)}\stackrel{\theta}{\longrightarrow}\dfrac{R\cap\gamma_{c+1}(F)}{\gamma_{c+1}(R,F,\dots,F)}.
\end{equation}
Moreover, since $M$ is $c$-central, so by part $(2)$ of Proposition \ref{pro1.4}, there is an epimorphism $\lambda:M\wedge^c\dfrac{L}{\gamma_{c+1}(L)}\longrightarrow T_{c+1}=\dfrac{\gamma_{c+1}(S,F,\dots,F)}{\gamma_{c+1}(R,F,\dots,F)}$. Hence 
\begin{equation}\label{equ7}
\dfrac{\gamma_{c+1}(S,F,\dots,F)}{\gamma_{c+1}(R,F,\dots,F)}\cong \dfrac{M\wedge^c\dfrac{L}{\gamma_{c+1}(L)}}{\ker\lambda}.
\end{equation}
Now, by putting \eqref{equ3}, \eqref{equ4}, \eqref{equ5}, \eqref{equ6}, and \eqref{equ7} in the sequence \eqref{Short-exact-seq.1}, we obtain the sequence \eqref{sequ.(b)}. 
\epr 
The following corollary is an immediate consequence of Proposition \ref{pro2.1} and helps us to prove the useful Proposition \ref{equalents}. 
\bc \label{cor2.2}
Let $L$ be a finite-dimensional $n$-Lie algebra and let $M$ be  its ideal. Then the following properties hold:
\begin{enumerate}
\item[(1).] 
The $n$-Lie algebra $\Mm^{(c)}(L)$ has dimension finite. 
\item[(2).]
$\dim\Mm^{(c)}(L/M)\leq \dim\Mm^{(c)}(L)+\dim\left(\dfrac{M\cap \gamma_{c+1}(L)}{\gamma_{c+1}(M,L,\dots,L)}\right)$.
\item[(3).]
$\dim\Mm^{(c)}(L)+\dim(M\cap \gamma_{c+1}(L))=\dim\Mm^{(c)}(L/M)+\dim\dfrac{\gamma_{c+1}(S,F,\dots,F)}{\gamma_{c+1}(R,F,\dots,F)}$,
where $F/R$ and $S/R$ are free presentations of $L$ and $M$, respectively. In particular, 
\[\dim\Mm^{(c)}(L)+\dim\gamma_{c+1}(L)=\dim\left(\gamma_{c+1}\left(\dfrac{F}{\gamma_{c+1}(R,F,\dots,F)}\right)\right).\]
\item[(4).] 
If $\Mm^{(c)}(L)=0$, then $\Mm^{(c)}(L/M)\cong \dfrac{M\cap\gamma_{c+1}(L)}{\gamma_{c+1}(M,L,\dots,L)}$.
\item[(5).]
If $M$ is an $r$-central ideal, that is, $M\subseteq Z_r(L)$ and $r\leq c$, then 
\[\dim\Mm^{(c)}(L)+\dim(M\cap \gamma_{c+1}(L))\leq\dim\Mm^{(c)}(L/M)+\dim\left(M\wedge^c\dfrac{L}{\gamma_{r+1}(L)}\right).\]
\end{enumerate}
\ec 
\bpr 
\begin{enumerate}
\item[(1).] 
It is clear according to the short exact sequence \eqref{Short-exact-seq.1} in the proof of Proposition \ref{pro2.1}.
\item[(2).] 
According to part $(1)$ of Proposition \ref{pro2.1}, this is clear.
\item[(3).] 
Since $\dfrac{\gamma_{c+1}(S,F,\dots,F)}{\gamma_{c+1}(R,F,\dots,F)}$ is an isomorphic image of $M\wedge^c\dfrac{L}{\gamma_{r+1}(L)}$, so by replacing $M\wedge^c\dfrac{L}{\gamma_{r+1}(L)}$ with $\dfrac{\gamma_{c+1}(S,F,\dots,F)}{\gamma_{c+1}(R,F,\dots,F)}$ in equation \eqref{sequ.(b)}, we have the following sequence:
\[0\longrightarrow \dfrac{\gamma_{c+1}(S,F,\dots,F)}{\gamma_{c+1}(R,F,\dots,F)}\longrightarrow\Mm^{(c)}(L)\longrightarrow\Mm^{(c)}(L/M)\longrightarrow M\cap \gamma_{c+1}(L)\longrightarrow 0.\]
Hence 
\[\dim\Mm^{(c)}(L)+\dim(M\cap \gamma_{c+1}(L))=\dim\Mm^{(c)}(L/M)+\dim\dfrac{\gamma_{c+1}(S,F,\dots,F)}{\gamma_{c+1}(R,F,\dots,F)}.\]
If $M=L$, then $\dim\Mm^{(c)}(L/M)=0$ and 
\[\gamma_{c+1}\left(\dfrac{F}{\gamma_{c+1}(R,F,\dots,F)}\right)=\dfrac{\gamma_{c+1}(F,F,\dots,F)}{\gamma_{c+1}(R,F,\dots,F)}.\] 
Therefore 
\[\dim\Mm^{(c)}(L)+\dim\gamma_{c+1}(L)=\dim\left(\gamma_{c+1}\left(\dfrac{F}{\gamma_{c+1}(R,F,\dots,F)}\right)\right).\]
\item[(4).] 
If $\Mm^{(c)}(L)=0$, then by part $(1)$ of Proposition \ref{pro2.1}, the proof is completed. 
\item[(5).] 
According to part $(2)$ of  Proposition \ref{pro2.1}, the result is obtained immediately.
\end{enumerate}
\epr 
The short exact sequence $0\longrightarrow M\longrightarrow L^*\longrightarrow L\longrightarrow 0$ of $n$-Lie algebras is called $c$-stem cover, if $M\subseteq Z_c(L^*)\cap \gamma_{c+1}(L^*)$ and $M\cong \Mm^{(c)}(L)$. Then $L^*$ is called $c$-cover of $L$. Specially, if $c=1$, then it is the definition of cover.   

Moreover, an $n$-Lie algebra $L$ is called $c$-perfect, if $\gamma_{c+1}(L)=L$.  
\bc \label{cor2.4}
Let $L$ be a $c$-perfect finite-dimensional $n$-Lie algebra with $c$-cover $L^*$. Then $\Mm^{(c)}(L)$ is trivial. 
\ec
\bpr 
By the definition of $c$-cover, there exists an ideal $M$ of $L$ such that 
\[L^*/M\cong L,\qquad M\subseteq Z_c(L^*)\cap \gamma_{c+1}(L^*),\qquad \Mm^{(c)}(L)\cong M.\]
Since $M$ is $c$-central and $L$ is $c$-perfect, that is, $L=\gamma_{c+1}(L)$, so by part (5) of \ref{cor2.2}, we have 
\[\dim\Mm^{(c)}(L^*)+\dim(M\cap \gamma_{c+1}(L^*))\leq \dim\Mm^{(c)}(L^*/M)+\underbrace{\dim(M\wedge^c \dfrac{L^*}{\gamma_{c+1}(L^*)})}_{0}.\]
Thus 
\[\dim\Mm^{(c)}(L^*)\leq \dim\Mm^{(c)}(L^*/M)=\dim\Mm^{(c)}(L).\]
Hence $\dim\Mm^{(c)}(L^*)\leq 0$, which implies that $\dim\Mm^{(c)}(L^*)=0$. 
\epr 
The authors in \cite[Proposition 1.6]{Akbarossadat-2-nilpotent} proved that if $L$ is an $n$-Lie algebra with the free representation $F/R$ and the map $\pi:F/\gamma_{c+1}(R,F,\dots,F)\longrightarrow F/R$ is an epimorphism, then $Z_c^*(L)=\pi\left(Z_c(F/\gamma_{c+1}(R,F,\dots,F))\right)$. Using this fact, we prove the following practical theorem.

\bp \label{equalents}
Let $L$ be a finite-dimensional $n$-Lie algebra and let $M$ be an ideal of $L$ with $M\subseteq Z_c(L)$. Then the following statements are equivalent:
\begin{enumerate}
\item 
$M\subseteq Z^*_c(L)$;
\item 
The natural map $\alpha:\Mm^{(c)}(L)\longrightarrow \Mm^{(c)}(L/M)$ is a monomorphism;
\item 
$\dim\Mm^{(c)}(L/M)=\dim\Mm^{(c)}(L)+\dim(M\cap\gamma_{c+1}(L))$.
\end{enumerate}
\ep 
\bpr 
Similar to \cite[Lemma 1.7]{Akbarossadat-2-nilpotent}, it can be proved that $M\subseteq Z^*_c(L)$ if and only if $\alpha:\Mm^{(c)}(L)\longrightarrow \Mm^{(c)}(L/M)$ is a monomorphism. 

If $M\subseteq Z^*_c(L)$, then by part (4) of Lemma \ref{lem2.1(2)} and part (2), we have the following sequence:
\[0\longrightarrow\Mm^{(c)}(L)\longrightarrow \Mm^{(c)}(L/M)\longrightarrow M\cap \gamma_{c+1}(L)\longrightarrow0.\] 
Hence 
\[\dim\Mm^{(c)}(L/M)=\dim\Mm^{(c)}(L)+\dim(M\cap\gamma_{c+1}(L)).\]

Conversely, if part (3) is hold, then by comparing part (3) of  Corollary \ref{cor2.2} and part (2), we get $\gamma_{c+1}(S,F,\dots,F)\subseteq \gamma_{c+1}(R,F,\dots,F)$. On the other hand, since $R\subseteq S$, so $\gamma_{c+1}(R,F,\dots,F)\subseteq \gamma_{c+1}(S,F,\dots,F)$. Therefore $\gamma_{c+1}(S,F,\dots,F)=\gamma_{c+1}(R,F,\dots,F)$, and this fact is equivalent to $S/\gamma_{c+1}(R,F,\dots,F)\subseteq Z_c(F/\gamma_{c+1}(R,F,\dots,F))$. From   \cite[Proposition 1.6]{Akbarossadat-2-nilpotent} we know that $Z_c^*(L)=\pi\left(Z_c(F/\gamma_{c+1}(R,F,\dots,F))\right)$, and hence \linebreak$\pi(S/\gamma_{c+1}(R,F,\dots,F))\subseteq Z_c^*(L)$. Therefore, $M\subseteq Z_c^*(L)$, which completes the proof. 
\epr 

The following proposition states that every $c$-perfect $n$-Lie algebra has at least one cover. 
\bp \label{pro2.3}
Let $L$ be a $c$-perfect $n$-Lie algebra. Then $L$ has at least one cover. 
\ep 
\bpr 
Let $\pi:F\longrightarrow L$ be an epimorphism with $\ker\pi=R$, and hence $F/R$ be a free presentation of $L$. Then 
\[0\longrightarrow\dfrac{R}{\gamma_{c+1}(R,F,\dots,F)}\longrightarrow\dfrac{F}{\gamma_{c+1}(R,F,\dots,F)}\stackrel{\bar{\pi}}{\longrightarrow}L\longrightarrow 0, \]
be the free $c$-central extension of $L$. 
Since $\bar{\pi}(\gamma_{c+1}(F)/\gamma_{c+1}(R,F,\dots,F))=L$, so by restricting $\pi$ to $\gamma_{c+1}(F)/\gamma_{c+1}(R,F,\dots,F)$, we have 
\begin{align*}
\ker\bar{\pi}=\ker\pi\bigcap \dfrac{\gamma_{c+1}(F)}{\gamma_{c+1}(R,F,\dots,F)} 
&=\dfrac{R}{\gamma_{c+1}(R,F,\dots,F)}\bigcap \dfrac{\gamma_{c+1}(F)}{\gamma_{c+1}(R,F,\dots,F)}\\
&=\dfrac{R\cap \gamma_{c+1}(F)}{\gamma_{c+1}(R,F,\dots,F)}=\Mm^{(c)}(L).
\end{align*}
Thus we get the following $c$-central extension:
\begin{equation}\label{equ-c-central-exten.}
0\longrightarrow\Mm^{(c)}(L)\longrightarrow \dfrac{\gamma_{c+1}(F)}{\gamma_{c+1}(R,F,\dots,F)}\longrightarrow L\longrightarrow 0.
\end{equation}
On the other hand, $L$ is $c$-perfect, so
\[L\cong \dfrac{F}{R}\longrightarrow F\cong L\oplus R\longrightarrow F=\gamma_{c+1}(F)+R.\]
Thus 
\begin{align*}
\gamma_{c+1}\left(\dfrac{\gamma_{c+1}(F)}{\gamma_{c+1}(R,F,\dots, F)}\right)&=
\dfrac{\gamma_{c+1}(\gamma_{c+1}(F),F,\dots, F)+\gamma_{c+1}(R,F,\dots, F)}{\gamma_{c+1}(R,F,\dots, F)}\\
&=\dfrac{\gamma_{c+1}(\gamma_{c+1}(F)+R,F,\dots, F)}{\gamma_{c+1}(R,F,\dots, F)}\\
&=\dfrac{\gamma_{c+1}(F)}{\gamma_{c+1}(R,F,\dots, F)}.
\end{align*}
Therefore, $\dfrac{\gamma_{c+1}(F)}{\gamma_{c+1}(R,F,\dots, F)}$ is $c$-perfect and so \eqref{equ-c-central-exten.} is a stem $c$-central extension of $L$.  
\epr 
In the following theorem, we give  upper and lower bounds for the dimension of $\Mm^{(c)}(L)$. 
\bt \label{the2.5}
Let $L$ be a finite-dimensional nilpotent $n$-Lie algebra of class $m$, let $d=d(L)$, and let $\dim\gamma_{i+1}(L)=a_i$, for all $i\geq 1$. Then 
the following inequalities are hold:
\[l_d^n(c+1)\leq\dim\Mm^{(c)}(L)+\dim\gamma_{c+1}(L)\leq l_d^n(c+1)+\sum_{i=1}^ca_id^{c(n-1)-i+1}.\]
If $L$ is abelian, then the lower and upper bounds are attained.
\et 
\bpr 
Since $L$ is nilpotent, so $\Phi(L)=L^2$, and hence $L/L^2$ is an abelian $n$-Lie algebra of dimension $d$. If we put $M=L^2$, then by Proposition \ref{pro1.2} we have 
\begin{equation}\label{equ-1}
\dim\Mm^{(c)}(L/M)=\dim\Mm^{(c)}(L/L^2)=l_d^n(c+1).
\end{equation} 
Thus by putting $M=L^2$ in part $(2)$ of Corollary \ref{cor2.2}, we obtain 
\begin{align}
\dim\Mm^{(c)}(L/L^2)  \leq  \dim\Mm^{(c)}(L)+\dim\dfrac{L^2\cap \gamma_{c+1}(L)}{\gamma_{c+1}(L^2,L,\dots, L)}
&=\dim\Mm^{(c)}(L)+\dim\dfrac{\gamma_{c+1}(L)}{\gamma_{c+2}(L)}\nonumber\\
&\leq \dim\Mm^{(c)}(L)+\dim\gamma_{c+1}(L) .\label{equ-2}
\end{align}
Therefore, by \eqref{equ-1} and \eqref{equ-2}, we have 
\begin{equation}\label{equ-(1)-a}
l_d^n(c+1)\leq\dim\Mm^{(c)}(L)+\dim\gamma_{c+1}(L).
\end{equation}
Now let $0\longrightarrow R\longrightarrow F\longrightarrow L\longrightarrow 0$ be a free presentation of $L$, and put $M=L^2$, in part $(1)$ of Proposition \ref{pro1.4}. Then 
\begin{equation*}
L=\dfrac{F}{R}\longrightarrow M=L^2=\left(\dfrac{F}{R}\right)'=\dfrac{F^2+R}{R},
\end{equation*}
and hence $T_{c+1}=\dfrac{\gamma_{c+1}(F^2+R,F,\dots,F)}{\gamma_{c+1}(R,F,\dots,F)}$, and it is an isomorphic image of $L^2\wedge^c L$, by part $(1)$ of Proposition \ref{pro1.4}. That is, there is an epimorphism $\psi:L^2\wedge^cL\longrightarrow T_{c+1}$ such that 
\begin{equation}\label{equ(3)}
\dfrac{L^2\wedge^cL}{\ker\psi}\cong T_{c+1}=\dfrac{\gamma_{c+1}(F^2+R,F,\dots,F)}{\gamma_{c+1}(R,F,\dots,F)}.
\end{equation}  
Now, if $M=L=F/R$, then 
$T'_{c+1}=\dfrac{\gamma_{c+1}(F)}{\gamma_{c+1}(R,F,\dots,F)}$. 
Thus 
\begin{align}
\dfrac{T'_{c+1}}{T_{c+1}}=\dfrac{\dfrac{\gamma_{c+1}(F)}{\gamma_{c+1}(R,F,\dots,F)}}{\dfrac{\gamma_{c+1}(F^2+R,F,\dots,F)}{\gamma_{c+1}(R,F,\dots,F)}}
&\cong \dfrac{\gamma_{c+1}(F)}{\gamma_{c+1}(F^2+R,F,\dots,F)}.\label{equ(*)}
\end{align}
On the other hand, by \eqref{equ(*)} we know that there exists the following epimorphism:
\begin{align*}
\theta: T'_{c+1}=\dfrac{\gamma_{c+1}(F)}{\gamma_{c+1}(R,F,\dots,F)} \stackrel{epi.}{\longrightarrow} \dfrac{T'_{c+1}}{T_{c+1}}=\dfrac{\gamma_{c+1}(F)}{\gamma_{c+1}(F^2+R,F,\dots,F)}.
\end{align*}
Then $\ker\theta=\dfrac{\gamma_{c+1}(F^2+R,F,\dots,F)}{\gamma_{c+1}(R,F,\dots,F)}$, and so we obtain the following short exact sequence:
\[0\longrightarrow \ker\theta\stackrel{inc.}{\longrightarrow}T'_{c+1}\stackrel{\theta}{\longrightarrow}\dfrac{T'_{c+1}}{T_{c+1}}\longrightarrow0.\]
That is, 
\begin{equation}\label{equ(4)}
0\longrightarrow \dfrac{\gamma_{c+1}(F^2+R,F,\dots,F)}{\gamma_{c+1}(R,F,\dots,F)}\stackrel{inc.}{\longrightarrow}\dfrac{\gamma_{c+1}(F)}{\gamma_{c+1}(R,F,\dots,F)}\stackrel{\theta}{\longrightarrow}\dfrac{\gamma_{c+1}(F)}{\gamma_{c+1}(F^2+R,F,\dots,F)}\longrightarrow0.
\end{equation}
Thus by \eqref{equ(3)} and \eqref{equ(4)}, we have 
\begin{equation}\label{equ(5)}
L^2\wedge^c L\longrightarrow\dfrac{\gamma_{c+1}(F)}{\gamma_{c+1}(R,F,\dots,F)}\stackrel{\theta}{\longrightarrow}\dfrac{\gamma_{c+1}(F)}{\gamma_{c+1}(F^2+R,F,\dots,F)}\longrightarrow0.
\end{equation}
Now, by part (3) of Corollary \ref{cor2.2}, and then using \eqref{equ(5)} and Theorem \ref{F^i/F^j}, we have 
\begin{align}
\dim\Mm^{(c)}(L)+\dim\gamma_{c+1}(L)&=\dim\left(\gamma_{c+1}\left(\dfrac{F}{\gamma_{c+1}(R,F,\dots,F)}\right)\right)\nonumber\\
&=\dim\left(\dfrac{\gamma_{c+1}(F)}{\gamma_{c+1}(R,F,\dots,F)}\right)\nonumber\\
&\leq \dim(L^2\wedge^c L)+\dim\left(\dfrac{\gamma_{c+1}(F)}{\gamma_{c+1}(\underbrace{F^2+R}_{\supseteq F^2},F,\dots,F)}\right)\nonumber\\
&\leq\dim(L^2\wedge^c L)+\dim\left(\dfrac{\gamma_{c+1}(F)}{\gamma_{c+1}(F^2,F,\dots,F)}\right)\nonumber\\
&=\dim(L^2\wedge^c L)+\dim\left(\dfrac{\gamma_{c+1}(F)}{\gamma_{c+2}(F)}\right). \label{equ(6)}
\end{align}
We know that $L^2\wedge^cL=(\dots((L^2\wedge L)\wedge L)\dots )\wedge L$ and also that the generators of  $L^2\wedge L$ are as $a_1\wedge l_2\wedge\dots \wedge l_n$, $a_1\wedge a_2\wedge l_3\wedge\dots\wedge l_n, \ldots,$ and $a_1\wedge a_2\wedge \dots\wedge a_{n-1}\wedge l_n$, where $a_i\in L^2$, and $l_j\in L$, for all $1\leq i\leq n-1$ and $2\leq j\leq n$. Indeed by using  equation \eqref{wedge3}, one can easily check that $a_1\wedge a_2\wedge l_3\wedge\dots\wedge l_n\in L^3\wedge L$, and similarly, $a_1\wedge \dots \wedge a_i\wedge l_{i+1}\wedge\dots\wedge l_n\in L^{i+1}\wedge L$. Thus 
\begin{equation}\label{equ(7)}
\dim(L^2\wedge^c L)\leq \sum_{i=1}^ca_id^{c(n-1)-i+1}.
\end{equation}  
Therefore, by \eqref{equ(6)} and \eqref{equ(7)}, we obtain 
\[\dim\Mm^{(c)}(L)+\dim\gamma_{c+1}(L)\leq l_d^n(c+1)+\sum_{i=1}^ca_id^{c(n-1)-i+1}.\]
\epr 
The next theorem states an upper bound for the dimension of $c$-nilpotent multiplier of finite-dimensional nilpotent $n$-Lie algebras of class $m$. 
\bt \label{the2.5-(2)}
Let $L$ be a $d$-dimensional nilpotent $n$-Lie algebra of class $m$. Then 
\begin{equation}\label{dimM(c)-m<c}
\dim\Mm^{(c)}(L)\leq \displaystyle\sum_{k=0}^{c}l_d^n(m+k),
\end{equation}
when $m\leq c$, and we have 
\begin{equation}\label{dimM(c)-m>c+1}
\dim\Mm^{(c)}(L)\leq \displaystyle\sum_{k=1}^{m}l_d^n(c+k),
\end{equation}
for all $c$ in which $m\geq c+1$, and if $L$ is abelian, then the equality is attained.
\et 
\bpr 
Let $L\cong \dfrac{F}{R}$. First assume that $m\leq c$. Since $\gamma_{m+1}(L)=0$, so 
\begin{align*}
0=\gamma_{m+1}(L)=\gamma_{m+1}\left(\dfrac{F}{R}\right)&=\dfrac{\gamma_{m+1}(F)+R}{R}\cong \dfrac{R}{\gamma_{m+1}(F)\cap R},
\end{align*}
which implies that  
\[\gamma_{m+1}(F)+R=\gamma_{m+1}(F)\cap R=R,\]
 that is ,
\[\gamma_{m+1}(F)\subseteq \gamma_{m+1}(F)\cap R\subseteq R.\]
Also, we know that $\gamma_{m+1}(F)\subseteq \gamma_{m}(F)$. Thus  
$\gamma_{m+1}(F)\subseteq \gamma_{m}(F)\cap R$. 
Therefore, 
\[\gamma_{c+m+1}(F)\subseteq \gamma_{c+m}(F)\subseteq \dots \subseteq \gamma_{m+2}(F)\subseteq \gamma_{m+1}(F)\subseteq R\cap \gamma_{m}(F).\] 
Hence
\begin{equation}\label{equ001}
\dim\dfrac{R\cap \gamma_m(F)}{\gamma_{c+m+1}(F)}\leq\dim\dfrac{\gamma_m(F)}{\gamma_{c+m+1}(F)}=\sum_{k=0}^{c}l_d^n(m+k).
\end{equation}
On the other hand, we know that 
\begin{align}
\dim\Mm^{(c)}(L)=\dim\dfrac{R\cap \gamma_{c+1}(F)}{\gamma_{c+1}(R,F,\dots,F)} &\leq \dim\dfrac{R\cap \gamma_m(F)}{\gamma_{c+1}(R,F,\dots,F)}\nonumber\\
&\leq \dim\dfrac{R\cap \gamma_m(F)}{\gamma_{c+1}(F)}\nonumber\\
&\leq \dim\dfrac{R\cap \gamma_m(F)}{\gamma_{c+m+1}(F)}. \label{equ002}
\end{align}
Thus by \eqref{equ001} and \eqref{equ002}, we obtain that 
\[\dim\Mm^{(c)}(L)\leq \sum_{k=0}^{c}l_d^n(m+k).\]
Now, suppose that $c+1\leq m$. Then 
\begin{align*}
\dim\Mm^{(c)}(L)=\dim\dfrac{R\cap\gamma_{c+1}(F)}{\gamma_{c+1}(R,F,\dots,F)}&=\dim\dfrac{\Big{(}R\cap \gamma_{c+1}(F)\Big{)}/\gamma_{m+1}(F)}{\gamma_{c+1}(R,F,\dots,F)/\gamma_{m+1}(F)}\\
&\leq \dim\dfrac{\gamma_{c+1}(F)/\gamma_{m+1}(F)}{\gamma_{c+1}(R,F,\dots,F)/\gamma_{m+1}(F)}\\
&\leq \dim\dfrac{\gamma_{c+1}(F)}{\gamma_{m+1}(F)}-\dim\dfrac{\gamma_{c+1}(R,F,\dots,F)}{\gamma_{m+1}(F)}\\
&\leq \dim\dfrac{\gamma_{c+1}(F)}{\gamma_{m+1}(F)}\\
&\leq \dim\dfrac{\gamma_{c+1}(F)}{\gamma_{c+m+1}(F)}=\sum_{k=1}^ml_d^n(c+k).
\end{align*}
\epr 
\bc \label{cor2.6}
Let $K$ be an $n$-Lie algebra over a field $\Bbb F$ of characteristic zero, and also suppose that the algebra $L=\dfrac{K}{Z_c(K)}$ satisfies in the stated conditions of Theorem \ref{the2.5}. Then 
\[\dim\gamma_{c+1}(K)\leq l_d^n(c+1)+\sum_{i=1}^ca_id^{c(n-1)-i+1}.\] 
\ec 
\bpr 
Let $F/R$ and $S/R$ be a free presentations of $K$ and $Z_c(K)$, respectively, where $S$ is an ideal of $F$. Therefore, according to Theorem \ref{the2.5}, we have
\begin{align}\label{equ01}
\dim\gamma_{c+1}(L)=\dim\gamma_{c+1}(K/Z_c(K))=\dim\gamma_{c+1}(F/S)=\dim\dfrac{\gamma_{c+1}(F)+S}{S}=\dim\dfrac{\gamma_{c+1}(F)}{\gamma_{c+1}(F)\cap S}.
\end{align}
Hence 
\begin{align}
\dim\gamma_{c+1}(L)=\dim\gamma_{c+1}(K/Z_c(K))&=\dim\dfrac{\gamma_{c+1}(K)+Z_c(K)}{Z_c(K)}\nonumber\\
&=\dim\gamma_{c+1}(K)+\dim Z_c(K)-\dim Z_c(K)=\dim\gamma_{c+1}(K).\label{equ02} 
\end{align}
Therefore, by equations \eqref{equ01} and \eqref{equ02}, we have 
\begin{align*}
\dim\gamma_{c+1}(K)=\dim\dfrac{\gamma_{c+1}(F)}{\gamma_{c+1}(F)\cap S}
&=\dim\left(\dfrac{\gamma_{c+1}(F)/\gamma_{c+1}(S,F,\dots, F)}{\gamma_{c+1}(F)\cap S/\gamma_{c+1}(S,F,\dots, F)}\right)\\
&\leq \dim\left(\dfrac{\gamma_{c+1}(F)}{\gamma_{c+1}(S,F,\dots, F)}\right)\\
&=\dim\Mm^{(c)}(K/Z_c(K))+\dim\gamma_{c+1}(K/Z_c(K))\\
&=\dim\Mm^{(c)}(L)+\dim\gamma_{c+1}(L)\\
&\leq  l_d^n(c+1)+\sum_{i=1}^ca_id^{c(n-1)-i+1}.
\end{align*}
\epr 
In the next proposition, we give another upper bound for the dimension of terms of lower central series, which is similar to   \cite[Lemma 14 ]{[13]}. 
\bp \label{pro2.7}
Let $L$ be an $n$-Lie algebra over a field $\Bbb F$ such that $\dim(L/Z_c(L))=d$. Then $\dim\gamma_{c+1}(L)\leq l_d^n(c+1)$. 
\ep 
\bpr 
Let $X=\{\bar{x}_1,\bar{x}_2,\dots,\bar{x}_d\}$ be a basis set for $L/Z_c(L)$. Then each member of $L$ can be shown as $\left(\sum\limits_{i=1}^dc_i\bar{x}_i\right)+z$, where $c_i\in\Bbb F$ and $z\in Z_c(L)$, for all $i=1,2,\dots, d$. If $F/R$ is a free presentation of $L$, then 
\begin{align*}
\dim\gamma_{c+1}(L)=\dim\gamma_{c+1}(F/R)=\dim\dfrac{F^{c+1}+R}{R}
=\dim\dfrac{F^{c+1}}{F^{c+1}\cap R}
&=\dim\dfrac{F^{c+1}/F^{c+2}}{(F^{c+1}\cap R)/F^{c+2}}\\
&\leq \dim\dfrac{F^{c+1}}{F^{c+2}}=l_d^n(c+1).
\end{align*} 
\epr 
The following corollary shows that, under some conditions, the converse of Proposition \ref{pro1.2} is hold. 
\bc \label{cor2.8}
Let $L$ be an $n$-Lie algebra with dimension $d$. Then the following properties hold:
\begin{enumerate}
\item[(1).] 
$\dim\Mm^{(c)}(L)+\dim\gamma_{c+1}(L)\leq l_d^n(c+1)$. 
In particular, 
\[\dim\Mm(L)\leq\dim\Mm(L)+\dim L^2\leq l_d^n(2).\] 
\item[(2).]
If $\dim\Mm^{(c)}(L)=l_d^n(c+1)$, then 
\begin{enumerate}
\item 
There exists an epimorphism $\theta:\Mm^{(c)}(L)\longrightarrow\Mm^{(c)}(L/L^2)$, which is obtained by choosing $M=L^2$ in the second part of Proposition \ref{pro2.1}. 
\item 
If $\ker\theta=0$, then $L$ is abelian. 
\end{enumerate}
\end{enumerate}
\ec 
\bpr 
\begin{enumerate}
\item[(1).] 
Put $K=F/\gamma_{c+1}(R,F,\dots,F)$, in  Proposition \ref{pro2.7}, where $F/R$ is a free presentation of $L$. Then 
\begin{align*}
\dim\dfrac{K}{Z_c(K)}&=\dim\left(\dfrac{F/\gamma_{c+1}(R,F,\dots, F)}{Z_c(F/\gamma_{c+1}(R,F,\dots, F))}\right)\\
&=\dim\left(\dfrac{F/\gamma_{c+1}(R,F,\dots, F)}{(Z_c(F)+\gamma_{c+1}(R,F,\dots, F))/\gamma_{c+1}(R,F,\dots, F)}\right)\\
&=\dim\left(\dfrac{F}{Z_c(F)+\gamma_{c+1}(R,F,\dots, F)}\right)\\
&\leq\dim \dfrac{F}{R}=d. 
\end{align*}
Therefore, 
\begin{equation}\label{equ1}
\dim\gamma_{c+1}(F/\gamma_{c+1}(R,F,\dots, F))\leq l_d^n(c+1).
\end{equation}
On the other hand, according to  part $(3)$ of Corollary \ref{cor2.2}, we have 
\begin{equation}\label{equ2}
\dim\Mm^{(c)}(L)+\dim\gamma_{c+1}(L)=\dim\gamma_{c+1}\left(\dfrac{F}{\gamma_{c+1}(R,F,\dots, F)}\right).
\end{equation}
Hence by \eqref{equ1} and \eqref{equ2}, we have
\[\dim\Mm^{(c)}(L)+\dim\gamma_{c+1}(L)\leq l_d^n(c+1).\]
In particular, if $c=1$, then 
$\dim\Mm(L)+\dim L^2\leq l_d^n(c+1)$, and so $\dim\Mm(L)\leq l_d^n(c+1)$.
\item[(2).] 
According to  Proposition \ref{pro2.1}, $\theta$ exists. Thus it is enough to show that $\theta$ is onto. 
By part $(1)$, we have 
\[\dim\Mm^{(c)}(L)+\dim\gamma_{c+1}(L)\leq l_d^n(c+1).\]
On the other hand, we know that $\dim\Mm^{(c)}(L)=l_d^n(c+1)$. So $\dim\gamma_{c+1}(L)=0$. Now, if we put $M=L^2$, in  part $(2)$ of  Proposition \ref{pro2.1}, then we obtain 
\[\dim (M\cap\gamma_{c+1}(L))=\dim (L^2\cap\gamma_{c+1}(L))=0,\]
and thus $\theta$ is an epimorphism. 

Now, let $\ker\theta=0$. Since $\theta$ is onto, so by putting $M=L^2$ in part $(2)$ of Proposition \ref{pro2.1}, we get the following exact sequence:
\[0\longrightarrow\Mm^{(c)}(L)\longrightarrow\Mm^{(c)}(L/L^2)\longrightarrow0,\]
that is, $\Mm^{(c)}(L)=\Mm^{(c)}(L/L^2)$. Since $L/L^2$ is abelian, so $L^2=0$ and hence $L$ is abelian.  
\end{enumerate}
\epr 

The following corollary determines the dimension of $c$-nilpotent multiplier of nilpotent $n$-Lie algebras of maximal class. We recall that a $d$-dimensional $n$-Lie algebra $L$ is nilpotent of maximal class $c+1$ when $\dim Z_c(L)=\dim \gamma_2(L)=d-n$ and also $\dim(Z_i(L)/Z_{i-1}(L))=1$. 
\bc \label{cor2.5(2)} 
Let $L$ be an $n$-Lie algebra of maximal nilpotency class $c+1$ with  dimension finite $d$.  Then 
\[\dim\Mm^{(c)}(L)\leq
 l_{d-1}^n(c+1)+n^{c(n-1)}.\]
\ec 
\bpr 
Since $L$ is an $n$-Lie algebra of maximal class $c+1$, 
then $Z(L)=\gamma_{c+1}(L)$, and hence $\dim\gamma_{c+1}(L)=1$. By using Theorem \ref{the2.4(2)}, if we choose $M=Z(L)$, then we obtain  
\begin{align*}
\dim\Mm^{(c)}(L)+\dim(\gamma_{c+1}(L)\cap Z(L))\leq &\dim\Mm^{(c)}(L/Z(L))+\dim\Mm^{(c)}(Z(L))\\
&+\dim((L/Z(L))^{ab}\, ^c\otimes Z(L)).
\end{align*}
On the other hand, we know that $Z_{c+1}(L)=L$ and that $\dim Z_c(L)=\dim\gamma_2(L)=d-n=c$. Moreover, for all $1\leq i\leq c$, we have $\dim Z_i(L)-\dim Z_{i-1}(L)=1$, and $\dim Z(L)=1$. Also, by  Proposition \ref{pro1.2}, $\dim\Mm^{(c)}(Z(L))=l_1^n(c+1)=0$, and also
$\dim(\gamma_{c+1}(L)\cap Z(L))=1$. Thus $Z(L)\subseteq \gamma_2(L)$. Hence,
\begin{align*}
\dim((L/Z(L))^{ab}=\dim\dfrac{L/Z(L)}{(L/Z(L))^2}&=\dim\dfrac{L/Z(L)}{(\gamma_2(L)+Z(L))/Z(L)}\\
&=\dim\dfrac{L/Z(L)}{\gamma_2(L)/Z(L)}\\
&=\dim (L/Z(L))-\dim(\gamma_2(L)/Z(L))\\
&=(d-1)-((d-n)-1)=n.
\end{align*}
Therefore,
\begin{align*}
\dim\Mm^{(c)}(L)+1&\leq \dim\Mm^{(c)}(L/Z(L))+0+n^{((c+1)-1)(n-1)}\times 1\\
&\leq  l_{d-1}^n(c+1)+n^{c(n-1)}.
\end{align*}
\epr 

Eshrati,   Saeedi, and  Darabi  \cite{Eshrati-Saeedi-Darabi} proved that every $d$-dimensional nilpotent $n$-Lie algebra of class $2$ with $\dim\gamma^2(L)=1$ has the structure $H(n,m)\oplus A(d-mn-1)$, where $H(n,m)$ is a Heisenberg $n$-Lie algebra of dimension $mn+1$ and $A(d-mn-1)$ is an abelian $n$-Lie algebra of dimension $d-mn-1$. Also, they gave the dimension of the Schur multiplier of $H(n,m)$ as follows:
\[\dim\Mm(H(n,m))=\begin{cases} 
n,&m=1,\\
{{mm}\choose{n}}-1,&m\geq2.
\end{cases}\]
Therefore, according to the above explanation and Proposition \ref{pro1.2}, if we get the dimension of the $c$-nilpotent multiplier of $H(n,m)$, then we will be able to calculate the dimension of the $c$-nilpotent multiplier of every finite-dimensional nilpotent $n$-Lie algebra of class $2$ with $\dim\gamma_2(L)=1$.
In the following theorem, we get the structure and dimension of the $c$-nilpotent multiplier of $H(n,m)$.
\bt \label{M(c)(H(n,m)}
If $L\cong H(n,m)$, then for all $c\geq2$, we have 
\[\Mm^{(c)}(H(n,m))\cong\begin{cases} 
A(l^n_{m}(c+1)+l^n_m(c+2)),&m=1,\\
A(l^n_{mn}(c+1)),&m\geq2.
\end{cases}\]
\et 
\bpr 
Assume first $m=1$, and let $F/R$ be a free presentation of $H(n,1)$. Since $cl(H(n,1))=2$, so $\gamma_{3}(H(n,1))=0$. Hence 
\[0=\gamma_3(H(n,1))=\gamma_3(F/R)=\dfrac{\gamma_3(F)+R}{R}\cong\dfrac{R}{\gamma_3(F)\cap R},\]
which  implies that $R=\gamma_3(F)$. Thus 
\[\Mm^{(c)}(H(n,1))=\dfrac{\gamma_{c+1}(F)\cap R}{\gamma_{c+1}(R,F,\dots,F)}=\dfrac{\gamma_{c+1}(F)\cap\gamma_3(F)}{\gamma_{c+1}(\gamma_3(F),F,\dots,F)}=\dfrac{\gamma_{c+1}(F)}{\gamma_{c+3}(F)}.\]
Therefore by  Theorem \ref{F^i/F^j}, we obtain  
\[\dim\Mm^{(c)}(H(n,1))=\dim\dfrac{\gamma_{c+1}(F)}{\gamma_{c+3}(F)}=l^n_{m}(c+1)+l^n_m(c+2).\]
Since the $c$-nilpotent multiplier of every $n$-Lie algebra is abelian, so $\Mm^{(c)}(H(n,1))\cong A(l^n_{m}(c+1)+l^n_m(c+2))$. 

Now, let $m\geq2$. Since $m\neq 1$, so $H(n,m)$ is not capable, by  \cite[Theorem 1.19]{Akbarossadat-2-nilpotent} and hence by  \cite[Lemma 1.5]{Akbarossadat-2-nilpotent}, we have $Z^*(H(n,m))\neq0$. Thus $Z^*(H(n,m))=\gamma_2(H(n,m))$, by   \cite[Proposition 1.10]{Akbarossadat-2-nilpotent}. Also, it is easy to check that $Z^*(H(n,m))\subseteq Z_c^*(H(n,m))$, and hence $\gamma_2(H(n,m))\subseteq Z_c^*(H(n,m))$. Thus by using Proposition \ref{equalents} and bearing in mind that $\gamma_{c+1}(H(n,m))=0$, so we have $\dim\Mm^{(c)}(H(n,m))=\dim\Mm^{(c)}(H(n,m))^{ab}$. Therefore, by Theorem \ref{pro1.2}, $\Mm^{(c)}(H(n,m))\cong\Mm^{(c)}(H(n,m))^{ab}\cong A(l_{mn}^n(c+1))$. 
\epr 
%

\ \\ \ \\

\noindent
{\small 
Farshid Saeedi$^*$: Corresponding author \\ 
Department of Mathematics, Mashhad Branch, Islamic Azad University, Mashhad, Iran.\\
E-mail address: saeedi@mshdiau.ac.ir \\ \ \\
Seyedeh Nafiseh Akbarossadat \\
E-mail address: n.akbarossadat@gmail.com  \\ 
Department of Mathematics, Mashhad Branch, Islamic Azad University, Mashhad, Iran.
}

\end{document}